\documentclass[12pt]{article}
\usepackage{yfonts}
\usepackage{eufrak}
\begin{document}
\def\F{\mathfrak{F}}
\def\L{\mathcal{L}}
\def\Rea{\mathrm{Re}\, }
\def\Ra{\mathfrak{R}}
\def\N{\mathbf{N}}
\def\C{\mathbf{C}}
\def\R{\mathbf{R}}
\def\bC{\mathbf{\overline{C}}}
\def\dist{\mathrm{dist}}
\def\cod{\mathrm{codim}\, }
\def\P{\mathbf{P}}
\title{On the Second Main Theorem of Cartan}
\author{Alexandre Eremenko\thanks{Supported by the NSF grant DMS-1361836.}}
\maketitle
\begin{abstract} The possibility of reversion
of the inequality in the 
Second Main Theorem of Cartan in the theory of holomorphic
curves in projective space is discussed. A new version of this theorem
is proved
that becomes an asymptotic equality for a class of
holomorphic curves defined by solutions of linear differential equations.

\end{abstract}

\section{Introduction}

We consider holomorphic curves $f:\C\to\P^n$. In homogeneous coordinates
such curves are represented as $(n+1)$-tuples of entire functions
$$f=(f_0:\ldots:f_n),$$ where not all $f_j$ are equal to $0$.
A homogeneous representation is called reduced if the $f_j$ do not have
zeros common to all of them. A reduced representation is defined up to
a common entire factor which is zero-free.

In the following definitions we use a reduced homogeneous representation,
however one can easily check that the definitions of
$N(r,a,f), T(r,f), N_1(r,f),$ $m(r,a,f)$ and
$m_k(r,f)$ are independent of the choice of a reduced
homogeneous representation.  

Let $a$ be a hyperplane in $\P^n$. It can be described by an
equation 
\begin{equation}
\label{hyperplane}
\alpha_0w_0+\ldots+\alpha_nw_n=0,\quad\mbox{where}\quad \alpha=(\alpha_0,\ldots,\alpha_n)\neq(0,\ldots,0).
\end{equation}
The intersection points of the curve $f(z)$
with the hyperplane $a$ are zeros of the entire function
$g_a=(\alpha,f)=\alpha_0f_0+\ldots+\alpha_nf_n$.
Let $n(r,a,f)$ be the number of these
zeros in the disc $|z|\leq r$, counting multiplicity,
then the Nevanlinna counting function is
defined as
\begin{equation}\label{N}
N(r,a,f)=\int_0^r(n(t,a,f)-n(0,a,f))\frac{dt}{t}+n(0,a,f)\log r.
\end{equation}
The Cartan--Nevanlinna characteristic $T(r,f)$ can be defined as follows:
\begin{eqnarray*}
T(r,f)&=&\frac{1}{2\pi}\int_0^r\left(
\int_{|z|\leq t}\Delta\log\| f(z)\| dm_z\right)\frac{dt}{t}\\
&=&\frac{1}{2\pi}\int_{-\pi}^\pi\log\| f(re^{i\theta})\|
d\theta-\log\| f(0)\|,
\end{eqnarray*}
where $\| f\|=\sqrt{|f_0|^2+\ldots+|f_n|^2}$,
and $dm$ is the element of the area.
Here $\Delta\log\| f\|$
is the density of the pull-back of the Fubini--Study metric, and equality
holds by Jensen's formula. The order $\rho$ of $f$ is defined by the formula
$$\rho=\limsup_{r\to\infty}\frac{\log T(r,f)}{\log r}.$$
The proximity functions are
$$m(r,a,f)=\frac{1}{2\pi}\int_{-\pi}^\pi
\log\frac{\| \alpha\|\| f(re^{it}\|}{|g_a(e^{it})|}dt.$$
Here the integrand is
$$\log\frac{1}{\dist(f(z),a)},$$
where $\dist$ is the ``chordal distance''
from the point $f(z)$ to the hyperplane $a$.
Now we consider the Wronskian determinant $W_f=W(f_0,\ldots,f_n)$ which
is an entire function; it is identically equal to zero if and only if 
$f$ is {\em linearly degenerate}, that is if
$f_0,\ldots,f_n$ are linearly dependent. We denote by $n_1(r,f)$
the number of zeros of $W_f$ in the disc $\{ z:|z|\leq r\}$ and define
the function $N_1(r,f)$ by a formula similar to (\ref{N}).

A set $A$ of hyperplanes is usually called admissible if
any $n+1$ hyperplanes
of the set have empty intersection. If the set $A$ contains
at least $n+1$ hyperplanes, admissibility is equivalent to
\begin{equation}\label{adm}
\cod (a_1\cap\ldots\cap a_k)=k
\end{equation}
for every $k\in[1,n+1]$ and every $k$ hyperplanes of the set $A$.
We use the convention that $\cod x=n+1$ iff $x=\emptyset$.
We use (\ref{adm}) to extend the definition of admissibility
to systems of arbitrary cardinality. So a system of hyperplanes
will be called
{\em admissible} if any $k\leq n+1$ vectors $\alpha$ defining these hyperplanes
as in (\ref{hyperplane}) are linearly independent.

With these definitions, the Second Main Theorem (SMT) of Cartan says: 
\vspace{.1in}

{\em For every linearly non-degenerate holomorphic curve
and for every finite admissible set $A$,
\begin{equation}\label{smt}
\sum_{a\in A}m(r,a,f)+N_1(r,f)\leq (n+1)T(r,f)+S(r,f),
\end{equation}
where $S$ is an ``error term'' with the property that
$S(r,f)=o(T(r,f))$ for $r\to\infty,\; r\not\in E$, where $E$ is an exceptional
set of finite length.}
\vspace{.1in}

Better estimates of the error term are available, but they do not concern us
here. When $n=1$, Cartan's SMT coincides with the Second Main Theorem
of Nevanlinna for the meromorphic function $f=f_1/f_0$. When $n=1$,
the assumption that
the set $A$ is admissible is vacuous.

Nevanlinna's SMT was considered from the very beginning as a partial
generalization
of the Riemann--Hurwitz formula \cite{Ahlfors}. However, the Riemann-Hurwitz
formula is an equality, while the SMT is only an inequality. This inspired
the research on the reversion of the SMT: roughly speaking,
the question is whether one
can replace the $\leq$ sign with the $=$ sign in (\ref{smt}) for $n=1$.
A survey of the early results on this topic is contained in the
book by Wittich \cite[Ch. IV]{Wittich}. The general conclusion one can make 
from these results is that for all simple, ``naturally arising'' meromorphic
functions an asymptotic equality indeed holds. But of course, (\ref{smt})
cannot be literally true for all meromorphic functions
in the form of equality, because there are
meromorphic functions $f$ with $m(r,a,f)\neq o(T(r,f))$ for an uncountable
set of $a$, and an exceptional set $E$ of $r$ does not help.

Recently, K. Yamanoi \cite{Y1} found a way to overcome this difficulty for
$n=1$. He defined the modified proximity function 
$$\overline{m}_q(r,f)=\sup_{(a_1,\ldots,a_q)\in\bC^q}\frac{1}{2\pi}\int_{-\pi}^\pi
\max_{1\leq j\leq q}\log\frac{1}{\dist(f(re^{it}),a_j)}dt.$$
With this definition, he proved the following theorem.
\vspace{.1in}

{\em Let $f:\C\to\bC$ be a transcendental meromorphic function. Let
$q:\R_{>0}\to\N$ be a function satisfying
$$q(r)\sim\left(\log^+\frac{T(r,f)}{\log r}\right)^{20}.$$
Then
$$\overline{m}_{q(r)}(r,f)+N_1(r,f)=2T(r,f)+o(T(r,f)),\quad r\not\in E,$$
where $E$ is a set of zero logarithmic density.}
\vspace{.1in}

For functions of finite order, this result was improved in \cite{Y2}: it
holds with any function $q(r)$ that satisfies $\log q(r)=o(T(r,f)).$

In this paper, we discuss the possibility of an asymptotic equality
in Cartan's SMT for arbitrary $n>1$. First we show by an example that the
admissibility condition creates a new difficulty which is not present
for $n=1$: even for very simple curves there can be no admissible
system for which (\ref{smt}) holds with equality.
Then we propose a modified form
of Cartan's SMT which does not
involve the admissibility condition, and show that in this modified form
asymptotic equality holds for a class of holomorphic curves.

\section{Example}

The simplest non-trivial examples in value distribution theory
for $n=1$ are meromorphic
functions $f=w_1/w_0$, where $w_0,\; w_1$ are
two linearly independent solutions of a differential equation of
the form
\begin{equation}\label{de}
w^{\prime\prime}+Pw=0,
\end{equation}
where $P$ is a polynomial. These functions $f$, which were
studied in detail
by F. Nevanlinna \cite{FN} and R. Nevanlinna \cite{RN}, can be characterized
by the properties: $f$ is of finite order, and $N_1(r,f)\equiv 0$.

For each such $f$, there is an integer $p$ and a finite set of
points $\{ a_1,\ldots,a_q\}$
in $\bC$ such that
\begin{equation}\label{777}
m(r,a_j,f)=(2m_j/p)T(r,f)+O(\log r),\quad r\to\infty,
\end{equation}
where $m_j$ are positive integers, and 
$$\sum_{j=1}^qm_j=p.$$
So we have an asymptotic equality in (\ref{smt}).

This result is related to two other results:

1. If $f$ is a meromorphic solution of arbitrary linear differential equation
with polynomial coefficients, then we have an asymptotic equality in
the SMT for $f$, with $A=\{0,\infty\}$, \cite[Ch. IV]{Wittich}.

2. If $f$ has finitely many critical and asymptotic values, then
an asymptotic equality holds in the SMT for $f$, if $A$ is the set
of critical and asymptotic values \cite{T,Wittich}. Functions $f=w_1/w_0$,
where $w_0,\; w_1$ are linearly independent solutions of (\ref{de})
have no critical values and their asymptotic values are exactly those $a_j$
in (\ref{777}). 

These results suggest
that in searching for improvements of (\ref{smt}) one has to
look first at the holomorphic curves whose homogeneous 
coordinates are linearly independent solutions of a differential equation
\begin{equation}
\label{den}
w^{(n+1)}+P_nw^{(n)}+\ldots+P_0w=0,
\end{equation}
with polynomial coefficients $P_j$. This class of curves can be 
characterized by the properties that the order is finite
and $N_1(r,f)\equiv 0$, \cite{P,E2}.

The following example was mentioned in \cite{E}:
\begin{equation}\label{de3}
w^{\prime\prime\prime}-zw^\prime-w=0.
\end{equation}
This is equivalent to
\begin{equation}
\label{airy}
w^{\prime\prime}-zw=c,\quad c\in\C.
\end{equation}
This is a non-homogeneous Airy equation, and we can describe the
asymptotic behavior of all solutions using the well-known asymptotic
formulas \cite{AS,W}. All non-trivial solutions are entire functions
of order $\rho=3/2$, and for description of their behavior we use the
Phragm\'en--Lindel\"of indicator:
$$h_w(t)=\lim_{r\to\infty}r^{-3/2}\log|w(re^{it})|.$$
First of all, we have three solutions $w_0,w_1,w_2$ (Airy's functions)
for $c=0$. These satisfy
\begin{equation}\label{0}
w_0+w_1+w_2=0,
\end{equation}
and have the indicators
\begin{equation}\label{8}
H_0(t)=-\cos\left(\frac{3}{2}t\right),\quad |t|\leq\pi,
\quad H_j(t)=H_0(t\pm2\pi/3),
\quad j=1,2.
\end{equation}
The rest of solutions of (\ref{de3}), which correspond to non-zero values
of $c$ in (\ref{airy}) can be expressed in terms of Airy functions
by the method of variation of constants. These explicit asymptotic expressions
show that the list of possible indicators for $c\neq 0$ is this:
\begin{equation}\label{9}
G_0(t)=\left(-\cos\left(\frac{3}{2}t\right)\right)^+,\quad |t|\leq\pi,\quad
G_j(t)=G_0(t\pm2\pi/3),\quad j=1,2.
\end{equation}
Another way to obtain these indicators is to notice that
(\ref{de3}) has a formal solution
$$w^*(z)=\sum_{n=0}^\infty\frac{(3n)!}{3^nn!}z^{-3n-1}.$$
According to the general theory \cite{W}, there exists a solution $w_3$
such that $w_3(z)$ has $w^*$ as the asymptotic expansion in
the sector
$$S_0=\{ z:|\arg z|<\pi/3\}.$$
For this solution, $h_{w_3}(t)=0,\, |t|\leq\pi/3$.
As the equation (\ref{de3}) is invariant under the
substitution $z\mapsto e^{2\pi i/3}z$,
in each of the three sectors $S_0,\, S_{\pm1}=e^{\pm2\pi i/3}S_0$
there exists a solution with zero indicator.

Notice that for every $t\not\in\{\pi,\pm\pi/3\}$,
the set of solutions with $h_w(t)\leq 0$ is at most
two dimensional. Indeed, if there were three linearly independent solutions
with $h_w(t)\leq 0$, then every solution would satisfy $h(t)\leq0$,
but this is
not so because $\max\{ H_0,H_1,H_2\}$ is positive at every point
$t\not\in\{\pi,\pm\pi/3\}.$
As for every $t$ there exists a solution $w$ with $h_w(t)=0$,
we obtain that for every $t$, the set of solutions $w$ with $h_w(t)<0$
is at most
one-dimensional. This shows that our list (\ref{8}), (\ref{9})
of possible indicators
of solutions is complete.

Now let $f$ be the holomorphic curve whose homogeneous coordinates
are three linearly independent solutions of (\ref{de3}). Then the
entire functions $g_a=(a,f)$
are exactly the non-trivial solutions of (\ref{de3}).
Let $A=\{ a_1,\ldots,a_q\}$ be an admissible system of hyperplanes.
Let $h_j$ be the indicators of entire functions $g_{a_j}$,
and let $h$ be their pointwise maximum. Then 
$$h(t)=|\cos((3/2)t)|,$$
$$\frac{1}{2\pi}\int_{-\pi}^\pi h(t)dt=
\frac{3}{2\pi}\int_{-\pi/3}^{\pi/3}\cos\left(\frac{3}{2}t\right)dt=
\frac{2}{\pi},$$
therefore
$$T(r,f)=\left(\frac{2}{\pi}+o(1)\right)r^{3/2}.$$
We claim that
$$\sum_{j=1}^q\int_{-\pi}^\pi(h(t)-h_j(t)dt\leq 8\int_{-\pi/3}^{\pi/3}
\cos\left(\frac{3}{2}t\right)=\frac{32}{3}.$$
This follows from the fact that on each of the three components of the set
$\{ t\in(-\pi,\pi):h(t)>0\}$ at most one of the $h_j$ can be negative,
and at most two of the $h_j$  can be non-positive, and in addition, we cannot
have negative indicators in all three components, because the three
solutions
$w_0,w_1,w_2$ satisfying (\ref{0}) cannot be all present in an admissible set.
So we have
$$\sum_{j=1}^qm(r,a_j,f)\leq\left(\frac{16}{3\pi}+o(1)\right)r^{3/2}\leq
\left(\frac{8}{3}+o(1)\right)T(r,f).$$
The Wronski determinant of three linearly independent solutions of
(\ref{de3}) is zero-free, $N_1(r,f)\equiv0$, and we cannot have asymptotic equality
in (\ref{smt}).
\vspace{.1in}

This example shows that if one desires (\ref{smt}) with asymptotic equality
then non-admissible sets of hyperplanes $A$ should be permitted.
In the next section we state and prove a version of (\ref{smt}) which applies
to an arbitrary finite system of hyperplanes.

\section{Modified Second Main Theorem}

Let us consider the projective space $\P^n$ equipped with
the chordal metric $\dist$. 
The distance between two subsets of
$\P^n$ is defined in the usual way, as the $\inf\dist(x,y)$, where
$x$ is in one set and $y$ is in another set.

Let us fix an {\em arbitrary} finite set $A$ of hyperplanes. Intersections
of various subsets of hyperplanes in $A$ are projective subspaces of various
codimension. We call these subspaces ``{\em flats} generated by $A$'', and
denote the set of all these flats by $F(A)$. We also denote by $\cod (x)$ the
codimension of a flat $x\in F(A)$. If $\cod(x)=k,$ then there exists
an admissible set $\{ a_1,\ldots,a_k\}\subset A$
such that $x=a_1\cap\ldots\cap a_k$.
If $\emptyset\in F(A)$, then flats of all codimensions 
$1,\ldots,n+1$ exist in $F(A)$.
Such systems $A$ will be called {\em complete}.
A system of hyperplanes is complete
if the vectors $\alpha$ corresponding to this system as in (\ref{hyperplane})
span $\C^{n+1}$.

We frequently use the following fact, without special mentioning:
if $a_1,\ldots,a_k$ is an admissible set of hyperplanes, and
$X=a_1\cap\ldots\cap a_k$ then
$$C_1\max_{1\leq j\leq k}\dist(w,a_j)\leq \dist(w,X)\leq C_2\max_{1\leq j\leq k}\dist(w,a_j),\quad w\in\P^n,$$
with positive
constants $C_1,C_2$ depending only on the set of hyperplanes.

For $w\in\P^n$ and $k\in\{ 1,\ldots,n\}$, we define $d_k(w)$
as the shortest distance from $w$ to a flat of codimension $k$ in $F(A)$.
It is also convenient to set $d_{n+1}(w)=1$.

For a holomorphic curve $f:\C\to\P^n$ and $k\in\{1,\ldots,n+1\}$,
we define the {\em $k$-proximity functions}
$$m_k(r,f)=\frac{1}{2\pi}\int_{-\pi}^\pi\log\frac{1}{d_k(f(re^{it}))}dt.$$
So $m_1\geq m_2\geq\ldots\geq m_n\geq m_{n+1}=0.$
Functions $m_k$ depend on $A$ which is not reflected in the notation.
Proximity functions for flats of arbitrary codimension were considered
for the first time by H. and J. Weyl's \cite{Weyls}.
With this definition we have
\vspace{.1in}

\noindent
{\bf Theorem 1.} {\em Let $f:\C\to\P^n$ be a linearly non-degenerate
holomorphic curve. 
Let $A$ be an arbitrary finite complete
set of hyperplanes. Then
\begin{equation}\label{msft}
\sum_{k=1}^nm_k(r,f)+N_1(r,f)\leq(n+1)T(r,f)+S(r,f),
\end{equation}
where $S(r,f)$ is the same error term as in Cartan's theorem.}
\vspace{.1in}

When $n=1$, we have 
$$m_1(r,f)=\frac{1}{2\pi}\int_{-\pi}^\pi\max_{a\in A}
\log\frac{1}{\dist(f(re^{it}),a)}dt+O(1),$$
so in the case when $m(r,a,f)=o(T(r,f))$ for all but finitely many $a$,
Theorem gives essentially the same as Yamanoi's result. 

Let $f_0,\ldots,f_n$ be linearly independent
 polynomials whose maximal degree is $k$.
Then there exist linear combinations $g_0,\ldots,g_n$ of these
polynomials whose degrees satisfy $k_0<k_1<\ldots<k_n=k$. 
Then we have
$$m_j(r,f)=(k-k_{j-1})\log r+O(1),\quad T(r,f)=k\log r+O(1).$$
Computing the degree of the Wronskian $W(g_0,\ldots,g_n)$, we obtain
\begin{equation}\label{wronskian}
N_1(r,f)=\left(\sum_{j=0}^nk_j-n(n+1)/2\right)\log r+O(1).
\end{equation}
Thus 
$$\sum_{j=1}^nm_j(r,f)+N_1(r,f)=(n+1)T(r,f)-\frac{n(n+1)}{2}\log r+O(1).$$
When $k$ is large, $T(r,f)$ is large in comparison with $\log r$,
and we obtain a relation close to (\ref{msft}).
So (\ref{msft}) can be considered as an extension of the formula for
the degree of the Wronskian to the transcendental case, compare
\cite[Introduction, ($\mathrm{II}^{\prime\prime}$)]{W}.

The proof of Theorem 1
is a combination of Cartan's argument with the following
elementary
\vspace{.1in}

\noindent
{\bf Lemma 1.} {\em Let $A$ be a finite complete set of hyperplanes
in $\P^n$. 
Then there exists a constant $C>0$ depending only on $A$, such that
for every $w\in\P^n$ we have
$$\prod_{k=1}^{n}d_k(w)\geq C\min_{B}\prod_{a\in B}\dist(w,a),$$
where the infimum is taken over all admissible systems
$B=\{ a_1,\ldots,a_{n+1}\}$
of hyperplanes in $A$.}
\vspace{.1in}

{\em Proof.} First we notice that if $x\in F(A)$ and $\cod x=k$,
then there exists an admissible subset $\{ a_1,\ldots,a_{n+1}\}\in A$
such that $x=a_1\cap\ldots\cap a_k$. Indeed, by passing from hyperplanes
to their defining vectors, this is equivalent to the familiar statement 
from linear algebra: if a finite set $A$ of vectors spans the space,
then every linearly independent subset of $A$ can be completed
to a basis consisting of vectors of $A$.

Now we prove the statement by contradiction.
For $w$ not in the union of hyperplanes of $A$, we set
$$\phi(w)=\frac{\prod_{k=1}^{n}d_k(w)}{\min_{B}\prod_{a\in B}\dist(w,a)}.$$
Suppose that there
is a sequence $w_j$ for which $\phi(w_j)\to 0$.
By choosing a subsequence, we may assume that $w_j\to w_\infty\in\P^n$.
If $w_\infty$ does not belong to any hyperplane $a\in A$, then
$\phi(w_\infty)>0$,
and we obtain a contradiction because $\phi$ is continuous in the 
complement of hyperplanes.

If $w_\infty$ belongs to some flat of $F(A)$,
let $x\in F(A)$ be the flat of maximal
codimension to which $w_\infty$ belongs. Then $d_j(w)$ are bounded away
from zero for $w$ in a neighborhood $V$ of $w_\infty$ and
$j>k=\cod x$. Suppose that $x=a_1\cap\ldots\cap a_k$. Then, by the remark
in the beginning, there is an admissible system
$B=\{ a_1,\ldots,a_{n+1}\}\subset A$ beginning with $a_1,\ldots,a_k$, and
$w_\infty\not\in a_j$ for $j>k$ by definition of $x$.
Then for $w\in V$, we have
$$\prod_{j=1}^nd_j(w)\geq C_1\prod_{j=1}^kd_j(w)\geq
C_2\prod_{j=1}^k\dist(w,a_j)
\geq C_3\prod_{j=1}^{n+1}\dist(w,a_j).$$
This contradicts
our assumption that $\phi(w_j)\to 0$ and proves the lemma.
\vspace{.1in}

{\em Proof of Theorem 1.} Fix a reduced representation of $f$. Normalize
all hyperplane coordinates so that $\| \alpha\|=1$ in (\ref{hyperplane}).
Let 
\begin{equation}
\label{df}
u=\log\| f\|,\quad u_a=\log|g_a|,\quad a\in A.
\end{equation}
Then 
$$-\log\dist(f(z),a)=u(z)-u_a(z).$$
According to Lemma 1, for every $z\in\C$, we can find
an admissible system $B(z),\, |B(z)|=n+1$, in $A$ such that
\begin{eqnarray}\nonumber
-\sum_{k=1}^n\log|d_k(f(z))|&\leq&-\sum_{a\in B(z)}\log\dist(f(z),a)+O(1)\\
&=&\label{one}
(n+1)u(z)-\sum_{a\in B(z)}u_a(z)+O(1).
\end{eqnarray}
Let $W=W(f_0,\ldots,f_n)$ be the Wronskian determinant.
If $B=\{ a_1,\ldots,a_{n+1}\}$ is an admissible system, then
\begin{equation}\label{w}
|W_B|=|W(g_{a_1},\ldots,g_{a_{n+1}})|=C(B)|W|.
\end{equation}
Let 
$$L_B(z)=\log^+\left|\frac{W_B(z)}{\prod_{a\in B(z)}|g_a(z)|}\right|.$$
Then
\begin{equation}\label{two}
-\sum_{a\in B(z)}u_a(z)\leq-\log|W_B(z)|+|L_B(z)|+O(1)\leq -\log|W(z)|+R(z),
\end{equation}
where $R(z)$ is the sum of non-negative quantities $L_B(z)$ over all
admissible systems of cardinality $n+1$.
The Lemma on the Logarithmic derivative
implies that
$$\int_{-\pi}^\pi R(re^{it})dt=S(r,f),$$
see \cite{C}, \cite[p. 222]{L}.
Jensen's formula gives 
$$\frac{1}{2\pi}\int_{-\pi}^\pi \log|W(re^{it})|dt=N_1(r,f)+O(1),$$
and the definition of $T(r,f)$ can be rewritten as
$$\frac{1}{2\pi}\int_{-\pi}^\pi u(re^{it})dt=T(r,f)+O(1).$$
Combining (\ref{one}) and (\ref{two}), integrating over circles
$|z|=r$, and using the last three equations
we obtain
$$\sum_{k=1}^nm_k(r,f)+N_1(r,f)\leq(n+1)T(r,f)+S(r,f).$$
This completes the proof of Theorem 1.
\vspace{.1in}

Now we compare Cartan's formulation of the SMT with Theorem 1.
\vspace{.1in}

\noindent
{\bf Proposition.} {\em Let $A$ be a finite admissible system of hyperplanes,
$|A|\geq n+1$, and $f$ a non-constant
holomorphic curve whose image
in not contained in any hyperplane of $A$.
Then
$$\sum_{a\in A}m(r,a,f)\leq\sum_{k=1}^nm_k(r,f)+O(1).$$
}
\vspace{.1in}

{\em Proof.} Let $A=\{ a_1,\ldots,a_q\}$. Define $u$ and $u_j=\log|g_{a_j}|$
by formulas (\ref{df}). Fix $z\in\C$ and order the functions $u_j$
by magnitude of $u_j(z)$,
$$u_{j_1}(z)\leq u_{j_2}(z)\leq\ldots\leq u_{j_q}(z),$$
where the $j_k$ depend on $z$.
Then for $k\leq n$ we have
$$u(z)-u_{j_k}(z)=-\log \dist(f(z),x_k)+O(1)\leq -\log d_k(f(z))+O(1),$$
where $x_k=a_{j_1}\cap\ldots\cap\ a_{j_k}$, 
and the $O(1)$ depends only on $A$. For $k\geq n+1$ we obtain
$u(z)-u_{j_k}(z)=O(1).$ Adding these
inequalities we obtain
$$\sum_{k=1}^q u(z)-u_{j_k}(z)\leq -\sum_{k=1}^n d_k(f(z))+O(1).$$
Integrating this inequality, over circles $|z|=r$
we obtain the statement of the proposition.
\vspace{.1in}

\noindent
{\bf Remark.} Unlike the usual proximity functions $m(r,a,f)$, the $m_k(r,f)$
can be substantially greater than $T(r,f)$. For example,
if $f$ is the curve considered in the previous section,
then $m_1(r,f)=2T(r,f)+O(1)$. It is a challenging problem to obtain
the exact upper estimates of the quantities
$$\delta_k=\liminf_{r\to\infty}\frac{m_k(r,f)}{T(r,f)}$$
for every $k\in[1,n]$. These are analogs of Nevanlinna defects. There is
a conjecture that $\delta_2\leq 1$ for $n=2$.

\section{Curves defined by solutions of linear ODE}

Let $\F$ be the set of all entire functions $y$ which satisfy 
differential equations of the form
\begin{equation}
\label{deN}
y^{(N)}+P_{N-1}y^{(N-1)}+\ldots+P_0y=0
\end{equation}
with polynomial coefficients $P_j$. 
This class contains 
exponential polynomials. For the curves of the form $f(z)=(e^{\lambda_0z}:
\ldots:e^{\lambda_nz})$ asymptotic equality holds in Cartan's SMT
\cite{A}.
\vspace{.1in}

\noindent
{\bf Theorem 2.} {\em Let $f:\C\to\P^n$ be a transcendental
linearly non-degenerate
holomorphic curve,
whose homogeneous
coordinates belong to $\F$. 
Then there exists a finite complete system $A$ of hyperplanes such that
$$\sum_{k=1}^nm_k(r,f)+N_1(r,f)=(n+1+o(1))T(r,f),\quad r\to\infty.$$
}
\vspace{.1in}

These curves are of finite order, so there is no exceptional set of $r$.
The result seems to be new even for $n=1$.

To prove Theorem 2, we use the following two facts
about the class $\F$:
\vspace{.1in}

\noindent
1. $\F$ is a differential ring \cite{FW}. This means that $\F$ is
closed under addition, multiplication and differentiation.

\vspace{.1in}

\noindent
2. For every differential equation (\ref{deN})
and every
$\theta$, there exists $\epsilon>0$, and $N$ linearly independent
solutions $y_1,\ldots,y_N$ of (\ref{deN}) such that
\begin{equation}\label{as}
y_k(z)\sim e^{Q_k(z^{1/p})}z^{s_k}\log^{m_k}z,\quad z=re^{it},\quad
r\to\infty,
\end{equation}
uniformly with respect to $t$ when $|t-\theta|\leq\epsilon.$
Here $Q_k$ are polynomials,
$Q_k(0)=0$, $p$ is a positive integer, $s_k\in\C$ and $m_k$ are
integers. All triples $(Q_k,n_k,m_k),\; 1\leq k\leq N,$ in (\ref{as})
are distinct. For a proof we refer to \cite{W}.
\vspace{.1in}

\noindent
3. All solutions $y$ of (\ref{deN}) are entire functions
of completely regular
growth in the sense of Levin--Pfl\"uger \cite{Levin},
the notion which we recall now.
\vspace{.1in}

Let $f$ be a holomorphic function in an angular sector
$S=\{ re^{i\theta}:|\theta-\theta_0|<\epsilon,\; r>0\}$. We say that $f$ has
{\em completely regular growth} with respect to order $\rho>0$
if the following finite limit exists
\begin{equation}\label{indic}
\lim_{r\to\infty,\; re^{i\theta}\not\in E}
\frac{\log|f(re^{i\theta})|}{|r|^\rho}=: h_f(\rho,\theta),
\end{equation}
uniformly with respect to $\theta$, for $|\theta-\theta_0|<\epsilon$.
Here $E\subset S$ is an exceptional set which can be covered by
discs centered at $z_j$ of radii $r_j$, such that
$$\sum_{j:|z_j|<r}r_j=o(r),
\quad r\to\infty.$$
Such sets $E$ are called $C_0$-sets in \cite{Levin}.

The limit $h_f(\rho,\theta)$ is called the indicator. It is always
continuous as a function of $\theta\in(-\epsilon,\epsilon)$.
Notice that if $f$ has completely regular growth with respect to
order $\rho$, then it has completely regular growth with respect
to any larger order, and the indicator with respect to the larger
order is zero.

An entire function $f$ is said to be of completely regular growth,
if it has completely regular growth in any sector with respect to its
order $\rho=\rho(f)$. 

If $f_1$ and $f_2$ are two functions of completely regular growth
with respect to the same order $\rho$ then evidently
$$h_{f_1+f_2}(\rho,\theta)\leq
\max\{ h_{f_1}(\rho,\theta),h_{f_2}(\rho,\theta)\},$$
and equality holds if $h_{f_1}(\rho,\theta)\neq h_{f_2}(\rho,\theta).$

Petrenko \cite[Sect. 4.3]{P} proved that all entire functions satisfying
differential equations of the form (\ref{deN}) have completely
regular growth.

Let $V\subset\F$ be a vector space of finite dimension $n+1$. Let $\rho$
be the maximal order of elements of $V$. From now on, all
indicators will be considered with respect to this
order $\rho$, and we
suppress it from notation.

Choose a ray
$\{ z:\arg z=\theta_0\}$. 
Each function $f\in V$ is a linear combination of some finite set
of entire functions $w_k$ which have asymptotics of the form
(\ref{as}) in an angular sector containing our ray. It is clear that
functions $w_k$ have trigonometric indicators of the form
$c_k\sin\rho(\theta-\theta_k)$. Two distinct trigonometric
functions of this form
can be equal only at a finite set of points.

We conclude that for each $V$ there exist finitely many rays such that
for any sector $S$ complementary to these rays the possible indicators
of elements of $V$ are strictly ordered:
\begin{equation}\label{sequence}
h_1(\theta)<h_2(\theta)<\ldots<h_m(\theta),\quad e^{i\theta}\in S.
\end{equation}
Here $m\geq 1$ is the number of distinct indicators in $S$.
Such sectors will be called {\em admissible} for a vector space $V$.

We fix an 
admissible sector $S$ of our vector space $V$,
and construct a {\em special basis}
in $V$. Let $h_j$ be the indicator of some element of $V$.
Then we define $V_j\subset V$ be the subspace consisting of
functions whose indicator at most $h_j$.
If all possible indicators are ordered as in (\ref{sequence}),
then $$V_1\subset V_2\subset\ldots\subset V_m=V.$$
We choose $\dim V_1$ linearly
independent functions in $V_1$, then $\dim V_2-\dim V_1$ functions
in $V_2$ which represent linearly independent
elements of the factor space $V_2/V_1$, and so on. So that the basis elements
chosen from $V_j\backslash V_{j-1}$ are linearly independent as
elements of $V_{j}/V_{j-1}$.

Let $w_0,w_1,\ldots,w_n$ be this basis, ordered in such a way that
the indicators increase,
\begin{equation}
\label{increase}
h_{w_0}(\theta)\leq h_{w_1}(\theta)\leq\ldots\leq h_{w_n}(\theta),\quad
e^{i\theta}\in S.
\end{equation}
Notice that, the indicator of any linear combination of the form
\begin{equation}\label{778}
c_0w_0+\ldots+c_{n-1}w_{n-1}+w_n
\end{equation}
is the same as the indicator of $w_n$. 
This sequence $(w_j)$ is called a special basis of $V$ corresponding
to the sector $S$.

\vspace{.1in}
\noindent
{\bf Lemma 2.} {\em Outside of a $C_0$
exceptional set $E$ as in (\ref{indic}), the special basis satisfies
$$\log|W(w_0,\ldots,w_n)|=\sum_{j=0}^n\log|w_j|+o(r^\rho)$$
in the sector $S$.}
\vspace{.1in}

{\em Proof}. If $f_1$ and $f_2$ are two functions of completely
regular growth in $S$, then the limit in (\ref{indic})
also exists for their ratio $f=f_1/f_2$ and this limit is equal
to $h_{f_1}(\theta)-h_{f_2}(\theta)$.
Let
$$\L(w_0,\ldots,w_n)=\frac{W(w_0,\ldots,w_n)}{w_0,\ldots,w_n}.$$
The statement of the Lemma is equivalent to $h_\L(\theta)\equiv 0.$

As $\L$ is a determinant consisting of the logarithmic derivatives
of functions of class $\F$, we have $h_\L(\theta)\leq 0$ by the
Lemma on the logarithmic derivative \cite{L}. It remains to prove
that $h_\L(\theta)\geq 0$.

We prove this by induction in $n$.
The statement is evident when $n=0$.
When $n=1$ we set $f=w_1/w_0$. Then $\L=f'/f$.
If $h_\L(\theta_0)<0$, we integrate $f'/f$ along the ray $\arg z=\theta_0$.
If the exceptional set $E$ intersects which ray, we bypass it by a curve
close to the ray consisting of arcs of circles. The result is that
$$f=c+O(e^{-\delta r^{\rho}}).$$
This implies that
$$h_{w_1-cw_0}(\theta_0)<h_{w_1}(\theta_0),$$
which contradicts the definition of the special basis.

Suppose now that the statement of the Lemma holds for spaces $V$
of dimension at most $m+1$, with some $m\geq 1$.
We have to prove it for $n=m+1$.
Assume by contradiction that
$h_{\L(w_0,\ldots,w_n)}(\theta_0)<0$ for some $\theta_0$. 
Define functions $B_j$ as solutions of the following system of linear equations
$$\sum_{j=0}^{n-1}B_jw^{(k)}_j=w^{(k)}_n,\quad k=0,\ldots,n-1.$$
By Cramer's rule,
$$B_j=\pm \frac{W_j}{W_n},$$
where $W_j$ is the Wronskian of size $n$ made of functions
$w_i$ with $i\neq j$.
We use the formula for differentiation of the logarithm of the
quotient of Wronskians
\cite[Part VII, Probl. 59]{PS}, \cite[p. 251]{L}
\begin{equation}\label{La}
\frac{d}{dz}\log\left(\frac{W_j}{W_n}\right)=\frac{W_{j,n}W}{W_jW_n}=
\frac{\L_{j,n}\L}{\L_j\L_n},
\end{equation}
where $W_{j,n}$ is the Wronskian of size $n-1$ with $w_j$ and $w_n$ deleted,
and $W$ is our Wronskian of size $n+1$. Notation $\L,\L_j,\L_{j,n}$
has similar meaning. Using the induction assumption,
we conclude that the right hand side of (\ref{La}) has negative
indicator.
Integrating with respect to $z$ along an appropriate curve
near the ray $\arg z=\theta_0$, that avoids the exceptional set $E$, we obtain
$B_j=c_j+O(e^{-\delta r^\rho}),\; 0\leq j\leq n-1$,
where $c_j\neq 0$ and $\delta>0$ are constants.
So we conclude that the indicator of
$$w_n-\sum_{j=0}^{n-1}c_jw_j$$
at the point $\theta_0$ is strictly less than $h_{w_n}(\theta_0)$.
This contradicts the property (\ref{778}) of the special basis.
The contradiction completes the proof of Lemma 2.
\vspace{.1in}

{\em Proof of Theorem 2.} Let $f:\C\to\P^n$ be a linearly
non-degenerate holomorphic curve whose
homogeneous coordinates are functions of $\F$.

Let $\rho$ be the order of our curve; it is equal
to the maximal
order of components $f_j$.

Let $V\subset \F$
be the subspace spanned by the homogeneous coordinates.
To such a space $V$ we associated finitely many exceptional rays,
whose complement consists of admissible sectors.
Let us fix any admissible sector $S$, and a special basis
$w_0,\ldots,w_n$ in $S$. 

Let $w_j=(f,\alpha_{j}), 0\leq j\leq n$, then the vectors
$\{ \alpha_0,\ldots,\alpha_n\}$ are
linearly independent. We define subspaces 
$$X_k=\{ w\in\C^{n+1}:(w,\alpha_0)=\ldots=(w,\alpha_{k-1})=0\},\quad
1\leq k\leq n,$$
so that $\cod X_k=k$.
We use the notation $u=\log\| f\|,\; u_j=\log|w_j|$.
If $z$ is outside of an exceptional set $E$,
we have
$$u_j(z)\leq u_{j+1}(z)+o(|z|^\rho),\quad 0\leq j\leq n-1,$$
view of (\ref{increase}). So
\begin{eqnarray*}
\log d_k(z)&\leq& \log\dist(f(z),X_k)\\
&=&
\max_{0\leq j\leq k-1}\log|(f(z),\alpha_j)|-\log \| f\|
=u_{k-1}(z)-u(z)+o(r^\rho).
\end{eqnarray*}
 Then, using Lemma 2 and $u=u_n+o(r^{\rho})$, we obtain
\begin{eqnarray*}
 \sum_{j=1}^{n} \log\frac{1}{d_k(z)} & \geq &
-\sum_{j=0}^{n-1}u_j(z)+nu+o(r^\rho)\\
&=&
-\sum_{j=0}^n u_j(z)+(n+1)u(z)+o(r^\rho)\\
&=&
-\log|W(w_0,\ldots,w_n)|+(n+1)u(z)+o(r^\rho).
\end{eqnarray*}
Integrating this with respect to $\theta$ on the sector $S$,
and then adding over
all admissible sectors, we obtain
$$\sum_{j=1}^n m_k(r,f)+N_1(r,f)\geq (n+1)T(r,f)+o(r^\rho).$$
Integrals over the exceptional set $E$ contribute $o(r^\rho)$ \cite{Levin}.
For curves $f$ with components in $\F$ we always have $T(r,f)=cr^\rho$,
so the error term is $o(T(r,f)$.

The opposite inequality follows from Theorem 1, where exceptional set
is absent because we deal with functions of finite order.
\vspace{.1in}

{\bf Remark.} A special case of Theorem 2 is that
the homogeneous coordinates of $f$ are linearly independent
solutions of (\ref{den}) with $N=n+1$. In this case we have $N_1(r,f)=0$.
For such curves Theorem
2 gives
$$\sum_{k=1}^nm_k(r,f)=(n+1+o(1))T(r,f).$$
These curves are analogous to meromorphic functions considered in
\cite{FN,RN}.
\vspace{.1in}

The author thanks Jim Langley for finding a mistake in the previous version
of this paper.

{\em Department of Mathematics

Purdue university

West Lafayette, IN 47907

eremenko@math.purdue.edu}
\end{document}